\documentclass[journal,transmag]{IEEEtran}

\usepackage{amsmath}
\usepackage{graphicx}
\usepackage{url}
\usepackage[margin=2.5cm,top=3.5cm,bottom=4cm,footskip=2cm,headsep=2cm]{geometry}
\usepackage{xcolor}

\def\div{\operatorname{div}}
\def\curl{\operatorname{curl}}
\def\Curl{\operatorname{\mathbf{Curl}}}

\def\B{\mathbf{B}}
\def\H{\mathbf{H}}
\def\J{\mathbf{J}}
\def\x{\mathbf{x}}
\def\y{\mathbf{y}}
\def\zero{\mathbf{0}}

\def\n{\mathbf{n}}

\title{On energy consistent vector hysteresis operators}

\author{Herbert Egger and Felix Engertsberger and Lukas Domenig and Klaus Roppert and Manfred Kaltenbacher}

\author{\IEEEauthorblockN{Herbert Egger\IEEEauthorrefmark{1,2}, Felix Engertsberger\IEEEauthorrefmark{1}, Lukas Domenig\IEEEauthorrefmark{3}, Klaus Roppert\IEEEauthorrefmark{3} and Manfred Kaltenbacher\IEEEauthorrefmark{3} 
}
\IEEEauthorblockA{\IEEEauthorrefmark{1}Institute for Numerical Mathematics, Johannes Kepler University, Linz, Austria}
\IEEEauthorblockA{\IEEEauthorrefmark{2}Johann Radon Institute for Computational and Applied Mathematics, Linz, Austria}
\IEEEauthorblockA{\IEEEauthorrefmark{3}Institute of Fundamentals and Theory in Electrical Engineering, TU Graz, Austria}
\thanks{Corresponding author: H. Egger (email: herbert.egger@jku.at).}}

\IEEEtitleabstractindextext{%
\begin{abstract}
Incremental models for magnetic vector hysteresis have been developed in previous works in accordance with basic principles of thermodynamics. In this paper, we present an equivalent representation of the associated hysteresis operator in terms of a co-energy functional which is useful for magnetic field computations based on a scalar potential. Using convex duality, we further define the corresponding energy functional and the associated inverse hysteresis operator which is required for computations based on the vector potential. The equivalence of the two representations with the energy-based hysteresis models proposed in earlier works is demonstrated and numerical results for some typical test problems  are presented obtained by finite element simulation of corresponding scalar and vector potential formulations.
\end{abstract}

}

\def\grad{\partial}

\begin{document}

\maketitle

\section{Introduction}
\label{sec:intro}

Constitutive models describing ferromagnetic  material behavior are typically based on an identity 
\begin{align} \label{eq:1}
\B = \mu_0 \H + \J
\end{align}
for magnetic induction $\B$, magnetic field intensity $\H$, and magnetic polarization $\J$, which is chosen here to represent the internal state of the system. An additional relation between $\J$ and either $\B$ or $\H$ is required to obtain models for particular classes of materials. 
In the presence of hysteresis~\cite{Bertotti1998,Mayergoyz2003}, a relation of the form 
\begin{align} \label{eq:2}
\J=\J(\H;\J_p)
\end{align} 
can be used to describe the change 
$\J_p \to \J$ in magnetic polarization from the previous to the current time step in response to a prescribed excitation field~$\H$. 
Such \emph{incremental} models of hysteresis thus implicitly establish a relation $\B=\B(\H;\J_p)$, which we call a \emph{forward hysteresis operator} in the following. This form of the material law is well-suited for the numerical solution of magnetic field problems based on the magnetic scalar potential formulations~\cite{BastosSadowski2014,Meunier2008}. 
For methods based on the vector potential formulation, on the other hand, the \emph{inverse hysteresis operator} $\H=\H(\B;\J_p)$ is formally required. Corresponding finite element implementations based on inverse Preisach and Jiles-Atherton models have been worked out in~\cite{Dupre1998,Gyselinck2000,Nierla2019} and \cite{Sadowski2002,LiKoh2015} for instance.

\subsection{Energy-based vector hysteresis models}
%
The hysteresis models of Bergqvist~\cite{Bergqvist1997,Bergqvist1997a}, later extended by Henrotte et al.~\cite{Henrotte2006,Henrotte2006b,Lavet2013} and Prigozhin et al.~\cite{Prigozhin2016}, provide particular definitions of the magnetic polarization $\J=\J(\H;\J_p)$ respecting basic principles of thermodynamics.
As shown in \cite{Kaltenbacher2022}, the above models are in fact equivalent.
The discretization of corresponding magnetic field problems has been briefly addressed in~\cite{Lavet2013,Prigozhin2016}, and fixed point methods have been used to solve the resulting nonlinear systems.
Faster Newton-type methods have been considered in~\cite{Jacques2016,Jacques2018}. The inverse relation $\H=\H(\B;\J_p)$ was obtained there by inverting the forward operator $\B=\B(\H;\J_p)$ numerically, which significantly complicates the implementation; see~\cite{Jacques2018} for an extensive discussion of this aspect. 
Since the energy-based hysteresis operators are not differentiable in a classical sense, the application of the standard Newton method seems questionable in general.
An efficient alternative, which can be justified rigorously, has been proposed recently in~\cite{Domenig2024,Egger2024quasi}. 

\subsection{Scope and main contributions}
%
We here provide an alternative characterization of the energy-based hysteresis operators $\B=\B(\H;\J_p)$ mentioned above in terms of a \emph{co-energy functional} that additionally depends on the magnetic polarization $\J_p$. 
This model is particularly useful for the analysis and numerical approximation of corresponding magnetic field problems based on the scalar potential formulation~\cite{Egger2024quasi,Engertsberger23}. 
Convex duality arguments allow us to define the associated \emph{energy functional} and to derive a compact form of the inverse hysteresis operator $\H=\H(\B;\J_p)$ which can be integrated naturally into discretization methods based on the vector potential. 
We further discuss the efficient implementation of the models and present some extensions which are required to improve their accuracy. 
Numerical tests are presented to illustrate the equivalence of the forward and inverse hysteresis operators and their performance in the context of finite element simulations.

\section{Forward hysteresis operator}
\label{sec:forward}

We start with recalling a simple instance of an energy-based hysteresis model describing ferromagnetic material behavior. Following the notation of~\cite{Lavet2013}, we define the polarization $\J=\J(\H;\J_p)$ as minimizer of the problem
\begin{align} \label{eq:3}
\min_\J \Big(U(\J) - \langle \H,\J\rangle + \chi |\J-\J_p| \Big).
\end{align}
Here $\J_p$ is the previous state of magnetic polarization, $\chi$ represents the strength of the pinning forces of the underlying microscopic description, and 
$U(\J)$ is an internal energy density, 
which is assumed to satisfy certain physically realistic structural conditions%
\footnote{$U(\J)$ is a proper convex function, continuously differentiable on its domain, and strongly convex. This in particular allows to verify existence of a unique minimizer.}; see~Section~\ref{sec:num} for a typical example. 
Together with the basic identity~\eqref{eq:1}, this leads to the material law
\begin{align} \label{eq:4}
\B 
= \mu_0 \H + \J(\H;\J_p),
\end{align}
with $\J(\H;\J_p)$ denoting the minimizer of~\eqref{eq:3}.  
the system~\eqref{eq:3}--\eqref{eq:4} is exactly the energy-based hysteresis model investigated in~\cite{Lavet2013}; see ~\cite{Bergqvist1997,Henrotte2006,Prigozhin2016} for equivalent formulations. 
The following key insight will be the basis for all our further considerations below.

\medskip 
\textbf{Assertion~1.}
\textit{The material law specified in~\eqref{eq:3}--\eqref{eq:4} can be stated equivalently in the form
\begin{align} \label{eq:5}
\B = \grad_\H w_*(\H;\J_p)
\end{align}
with co-energy density $w_*(\H;\J_p)$ defined implicitly by
\begin{align} \label{eq:6}
&w_*(\H;\J_p) \\
&=\frac{\mu_0}{2} |\H|^2 
- \min_\J\Big( U(\J) 
- \langle \H,\J\rangle + \chi |\J-\J_p| \Big). \notag
\end{align}}

\smallskip
This claim follows by application of Danskin's theorem%
\footnote{Let $f(\x,\y)$ be a sufficiently smooth functions and strongly convex in $\y$. Then the function $\phi(\x) = \inf_\y f(\x,\y)$ is well-defined, i.e., for every $x$ there exists a unique minimizer $\y(\x)$ of $f(\x,\y)$, and $\phi(\x)$ is differentiable with $\partial_\x \phi(\x) = \partial_\x f(\x,\y(\x))$. If $f$ is continuously differentiable in both variables, this is clear from the implicit function theorem. Let us note that the regularity of $f(\x,\y)$ with respect to $\y$ can be further relaxed substantially~\cite{Danskin1966}.}
and some elementary computations, which yield
\begin{align} \label{eq:7}
\grad_\H w_*(\H;\J_p) = \mu_0 \H + \J(\H;\J_p)
\end{align}
with $\J(\H;\J_p)$ denoting the minimizer of~\eqref{eq:3}. 
\medskip 
\textit{Remark~1.}
The material law~\eqref{eq:5} is well-suited for the scalar potential formulation of magnetic field problems which requires a material law of the form $\B=\B(\H)$. The availability of the co-energy density~$w_*$ is also useful for the analysis and numerical solution of the corresponding nonlinear partial differential equations~\cite{Engertsberger23}.

\medskip 

For application in a vector potential formulation, access to the \emph{inverse hysteresis operator} $\H=\H(\B)$ is required instead; see e.g.~\cite{Jacques2016}.
In the following we derive a compact form of this inverse relation for~\eqref{eq:5}, which opens the door for a systematic analysis and and efficient numerical solution of corresponding magnetic field problems based on the vector potential formulation. 

\section{Inverse hysteresis operator}
\label{sec:inverse}

From general considerations in material modeling and convex analysis~\cite{Silvester1991,BorweinLewis2000}, we know that the inverse of the constitutive law \eqref{eq:5} can be stated in the form
\begin{align} \label{eq:8}
\H = \grad_\B w(\B;\J_p),
\end{align}
where $w(\B;\J_p)$ denotes the energy density, i.e., the convex-conjugate of the co-energy density $w_*$ with respect to the first argument. This function is defined by
\begin{align} \label{eq:9}
w(\B;\J_p) = \sup\nolimits_\H \langle \H,\B\rangle - w_*(\H;\J_p).   
\end{align}
Recall that $w_*(\H;\J_p)$ was already defined implicitly via a minimization problem, hence this formula is mainly useful for theoretical purposes. 
We now derive an equivalent representation for the energy density~$w(\B;\J_p)$ which is better suited for analysis and implementation.

\medskip 
\textbf{Assertion~2.}
\textit{The energy density~\eqref{eq:9} can be written as
\begin{align} \label{eq:10}
&w(\B;\J_p) \\
&= \min_\J \Big(\frac{1}{2\mu_0} |\B-\J|^2 + U(\J) + \chi |\J-\J_p| \Big). \notag
\end{align}}

By application of Danskin's theorem, we here obtain
\begin{align} \label{eq:11}
\grad_\B w(\B;\J_p) = \frac{1}{\mu_0} (\B-\tilde \J(\B;\J_p)) 
\end{align}
where $\tilde \J(\B;\J_p)$ denotes the unique minimizer in~\eqref{eq:10}. Using~\eqref{eq:1} and comparing with~\eqref{eq:4} then readily shows that $\tilde \J(\B;\J_p) = \J(\H;\J_p)$ and hence validity of~\eqref{eq:9}.

\medskip 
\textit{Remark~2.}
The evaluation of $w_*(\H;\J_p)$  and $w(\B;\J_p)$ each require the solution of a convex minimization problem of similar form. 
After determining the minimizers $\J=\J(\H;\J_p)$ resp. $\J=\tilde \J(\B;\J_p)$, the computation of $\B=\grad_\H w_*(\H;\J_p)$ and $\H=\grad_\B w(\B;\J_p)$  consists of a simple evaluation of the formulas \eqref{eq:5} or \eqref{eq:11}.
These are the key steps for the implementation of the constitutive models required during simulation.
The numerical complexity of the two models, i.e., the forward and inverse hysteresis operator, therefore is essentially the same. 

\section{Multiple pinning forces}
\label{sec:multiple}

As noted in \cite{Bergqvist1997,Henrotte2006,Prigozhin2016}, the accuracy of the respective energy-based vector hysteresis models can be significantly increased by considering multiple pinning forces and, correspondingly, multiple internal variables. 
We now discuss these generalizations in our framework.

\subsection{Forward hysteresis operator}
%
Following the notation of~\cite{Jacques2016}, we split 
\begin{align} \label{eq:12}
\J = \sum\nolimits_k \J_k
\qquad \text{and} \qquad 
\J_p = \sum\nolimits_k \J_{p,k}
\end{align}
into partial polarizations $\J_k$ and $\J_{p,k}$, respectively. 
We further use $\{\J\}$ and $\{\J_{p}\}$ to denote the  collections of the respective components. 
In accordance with the previous formula~\eqref{eq:6}, we now define the co-energy density by
\begin{align} \label{eq:13}
&w_*(\H;\{\J_{p}\}) = \frac{\mu_0}{2} |\H|^2 \\
& -\sum\nolimits_k \min_{\J_k} \Big( U_k(\J_k) - \langle \H, \J_k\rangle + \chi_k |\J_k - \J_{p,k}| \Big). \notag
\end{align}
With the very same reasoning as before, we obtain the following generalization of our first result.

\medskip 
\textbf{Assertion~3.}
\textit{The material law~$\B=\partial_H w_*(\H;\{\J_p\})$, compare with~\eqref{eq:5} for one pinning force, can be stated as 
\begin{align} \label{eq:14}
\B = \mu_0 \H + \sum\nolimits_k \J_k,
\end{align}
where $\J_k=\J_k(\H;\J_{p,k})$ are the minimizers in~\eqref{eq:13}}. 

\medskip 
Let us emphasize that this model is again equivalent to the generalization of the vector-hysteresis models considered in~\cite{Lavet2013,Jacques2016} with  multiple pinning forces. 

\subsection{Inverse hysteresis operator}
%
By generalizing~\eqref{eq:10}, we obtain the corresponding energy density $w(\B;\{\J_p\})$ for multiple pinning forces
\begin{align} \label{eq:15}
&w(\B;\{\J_p\}) 
= \min_{\J_k} \Big( \frac{1}{2\mu_0} \big|\B - \sum\nolimits_k \J_k\big|^2  \\
&\qquad \qquad  + \sum\nolimits_k U_k(\J_k) + \chi_k \big|\J_k - \J_{p,k}\big| \Big). \notag
\end{align}
With similar reasoning as in our previous considerations, we thus immediately obtain our final theoretical result.

\medskip 
\textbf{Assertion~4.}
\textit{The material law $\H = \partial_\B w(\B;\{J_p\})$, see~\eqref{eq:8} for one pinning force, can be written as 
\begin{align} \label{eq:16}
\H=\frac{1}{\mu_0} (\B - \sum\nolimits_k \J_k)
\end{align}
where  $\J_k=\tilde \J_k(\H;\{\J_p\})$ are the minimizers in~\eqref{eq:15}. 
Moreover, the function $w(\B;\{\J_p\})$ stated in~\eqref{eq:15} is related to $w_*(\H;\{\J_p\})$ given in~\eqref{eq:13} by convex duality~\eqref{eq:9}.} 

\medskip
The claims follow by comparing~\eqref{eq:14} and \eqref{eq:16}, and using the very same arguments as for the analysis of the inverse hysteresis operator with a single pinning force. 

\subsection{Computational complexity}
%
The partial polarization $\J_{k}=\J_k(\H;\J_{p,k})$ in the co-energy functional~\eqref{eq:13} can be computed efficiently by solving independent optimization problems for each~$k$.
The determination of $\J_k=\tilde \J(\B;\{\J_p\})$ in the energy functional~\eqref{eq:15}, on the other hand, does not separate into independent minimization problems for each $k$. 
For multiple pinning forces, the computational effort for evaluating the inverse law~\eqref{eq:16} thus is somewhat larger than that for computing the forward relation~\eqref{eq:14}; see Table~\ref{table:1} below for computational results.
Up to possible round-off errors, both variants, however, still represent the same material law. 
Further note that a similar coupling between the inner minimization problems for determining the partial polarizations $\J_k$ occurs already also in the forward hysteresis operator, if an effective field $\H_\text{eff} = \H + \alpha \J$ is used in the formulation of the model; we refer to~\cite{Prigozhin2016} for comments in this direction.

\section{Numerical validation}
\label{sec:num}
We now illustrate the validity of Assertions~1--4 by some numerical results for the forward and inverse hysteresis models with one and multiple pinning forces. 

\subsection{General considerations}
%
The evaluation of the magnetic co-energy density~\eqref{eq:13} and the magnetic energy density~\eqref{eq:15}, both, require the solution of convex non-smooth minimization problems.
To avoid technicalities and to provide a unified numerical approach to all problems under consideration, we regularize the non-smooth term $|\J - \J_{p}| \approx |\J - \J_{p}|_\epsilon$, where $|\x|^2_\epsilon := |\x|^2 + \epsilon$ with $\epsilon = 10^{-8}$ chosen in our tests. 
The resulting problems are then treated by a Newton-method with line search~\cite{Nocedal2006}. 
Following~\cite{Prigozhin2016}, specialized algorithms for minimizing the original non-smooth function could be devised when the minimization problem separates into independent smaller minimization problems. 
To allow for a better comparison of the different models, we refrain from using such specialized solution approaches in our numerical tests. 
All computations were carried out in \textsc{Matlab} and on a desktop computer with 13th Gen Intel Core i5-1335U CPU with clock rate 1.30 GHz.
\subsection{Internal energy densities}
%
The choice of our test problems is inspired by the ones in \cite{Prigozhin2016}. Following this reference, we choose 
\begin{align} \label{eq:Uk}
U_k(\J) = - \frac{A w_k J_s}{\pi} \log \left(\cos \left(\frac{\pi}{2} \frac{|\J|_\epsilon}{w_k J_s}\right)\right)
\end{align}
for the internal magnetic energy densities.
Here $J_s$ is the saturation polarization, $A$ corresponds to the steepness of the anhysteretic curve, and $w_k$ is a weight factor for the pinning forces $\chi_k$.
Specific values for $A$, $J_s$ and $w_k$, $\chi_k$ will be provided for each test.
In the case of one single pinning force, we drop the index $k$.

\subsection{Single material point (one pinning force)}
\begin{figure}[ht!]
    \centering
    \includegraphics[trim={1cm 0cm 1cm 0cm},clip,width=0.48\textwidth]{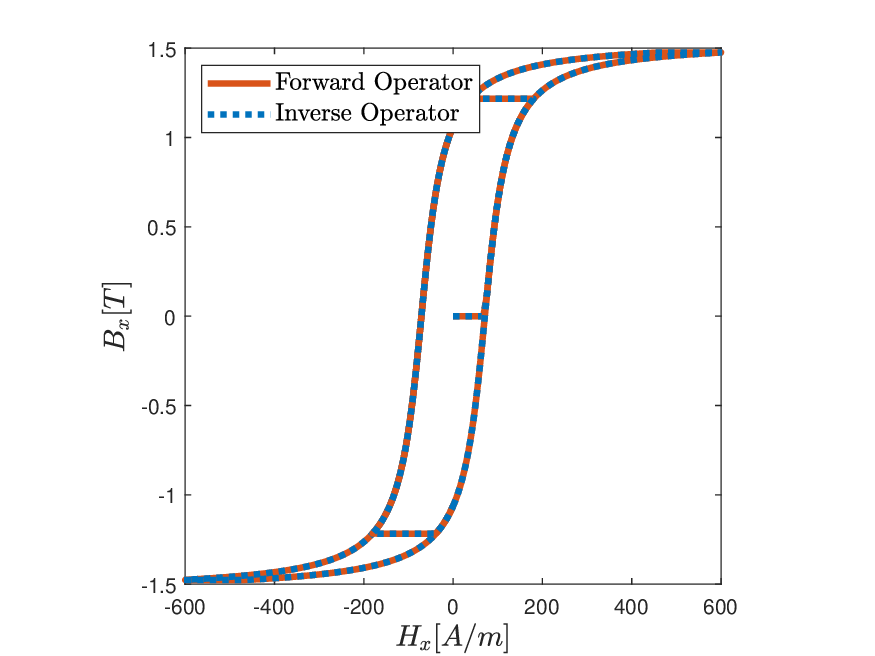}
    \caption{Hysteresis loop for one pinning force $(\chi_1 = 71, w_1 = 1)$ with $\H^i = (H_m\,\text{sin}(t^i),0)$, where $H_m = \{180,600\}$.}
    \label{fig:1}
\end{figure}

\begin{figure}[ht!]
    \centering
    \includegraphics[trim={1cm 0cm 1cm 0cm},clip,width=0.48\textwidth]{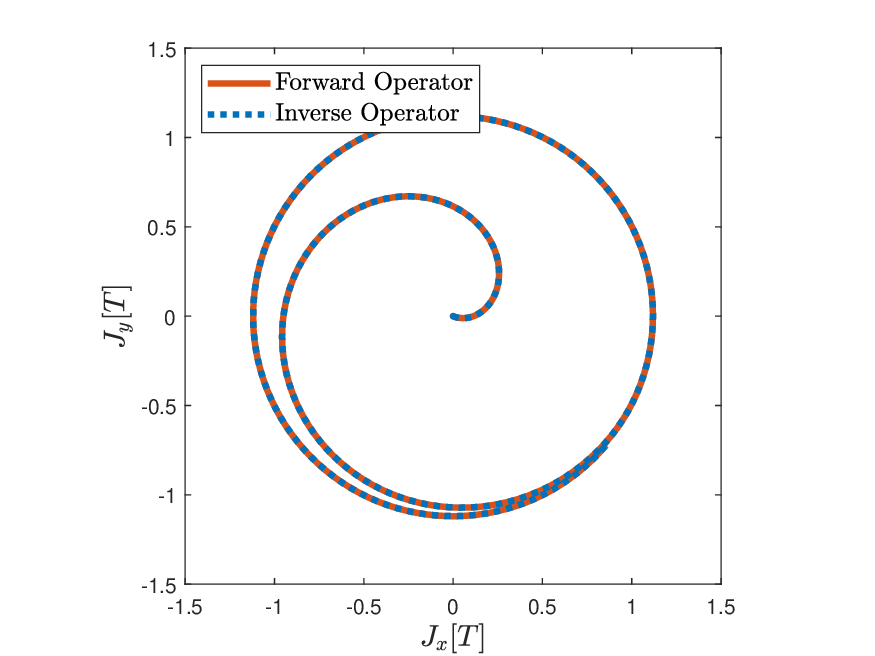}
    \caption{Hysteresis loop for one pinning force $(\chi_1 = 71, w_1 = 1)$ with $\H^i = (H_m^i\,\text{sin}(t^i),\text{cos}(t^i))$, where $H_m^i = 110\min(t^i/6\pi,1)$.}
    \label{fig:2}
\end{figure}

\noindent
We start with reproducing the results of \cite[Fig.~1,2]{Prigozhin2016}.
The parameters $A = 38, J_s = 1.54$\,T and $\chi = 71$, $w=1$ are used in the definition \eqref{eq:Uk} of the internal energy density $U=U_k$. 
A series of magnetic field vectors $\H^i$, $i=1,2,\ldots$ is defined, which serve as input for the forward hysteresis operator~\eqref{eq:5}. 
The magnetic polarization $\J^i=\J(\H^i;\J_{p}^i)$ are determined by solving the inner minimization problems~\eqref{eq:4}, with $\J_p^i=\J^{i-1}$ and $\J^0=\zero$ is chosen for initialization.
The output of this model are magnetic flux vectors $\B^i$.  
The red curves in Fig.~\ref{fig:1} and \ref{fig:2} display the results of these computations, corresponding to the results in \cite[Fig.~1,2]{Prigozhin2016}. 

In a second test run, we use the output vectors $\B^i$ of the previous simulation as input vectors for the inverse hysteresis operator~\eqref{eq:9}. The magnetic polarizations $\J^i=\tilde \J(\B^i;\J_p^i)$ are determined by solving the inner minimization problems in \eqref{eq:11} with $\J_p^i=\J^{i-1}$ and $\J^0=\zero$ as before.
The vectors $\B^i$ are plotted against the output vectors $\H^i$ of these computations by blue lines in Fig.~\ref{fig:1}~and~\ref{fig:2}.
In accordance with Assertion~1~and~2, excellent agreement with the results of \cite{Prigozhin2016} and a perfect match of the results obtained by the forward and inverse hysteresis models \eqref{eq:5} and \eqref{eq:8} is obtained. 

\subsection{Single material point (multiple pinning forces)}
\begin{figure}[ht!]
    \centering
    \includegraphics[trim={1cm 0cm 1cm 0cm},clip,width=0.45\textwidth]{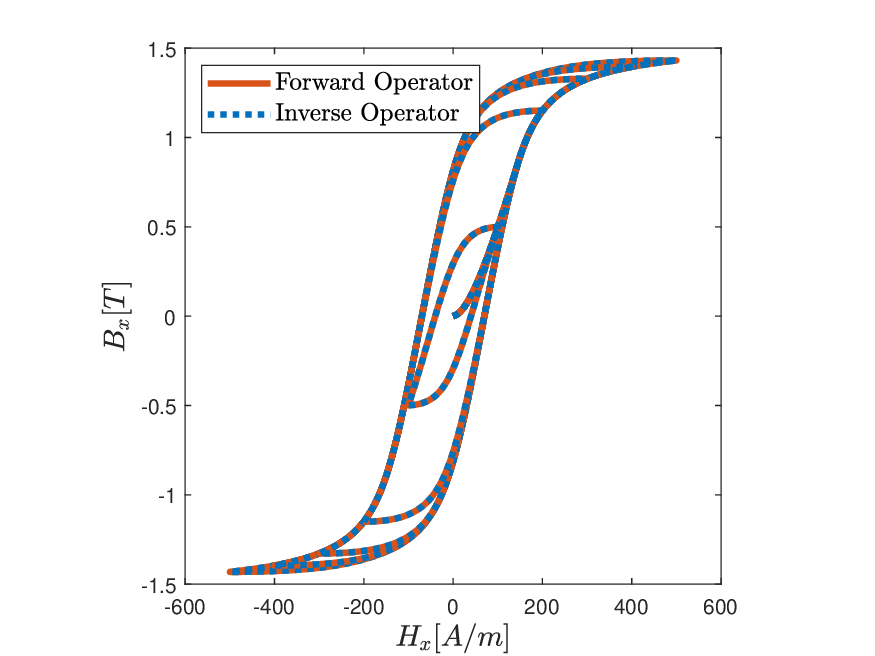}
    \caption{Hysteresis loop for multiple pinning forces ($N_{\chi} = 20$) with $\H^i = (H_m\,\text{sin}(t^i),0)$, where $H_m = \{100,200,300,400,500\}$.}
    \label{fig:3}
\end{figure}

As a second example, we consider a more realistic composite model involving multiple pinning forces. 
For comparison with previous work, we consider the setting of \cite[Fig.~4]{Prigozhin2016}. The internal energy densities \eqref{eq:Uk} are used with the parameters $A = 50$, $J_s = 1.54$\,T, and $\chi_k = 140(k-1)/(N_{\chi}-1)$, $w_k=1/20$ for $k = 1,\dots N_{\chi}$.
The red lines in Fig.~\ref{fig:3} display the results obtained by application of the forward hysteresis operator \eqref{eq:14} with given vectors $\H^i$ specified as inputs.  
The corresponding results of the inverse hysteresis operator~\eqref{eq:16}, which uses the output vectors $\B^i$ of the forward model as inputs, are displayed in blue. 
Again the results obtained the forward and inverse hysteresis operator align perfectly, 
clearly demonstrating the equivalence of the two models. In addition, the results match well with those reported in~\cite{Prigozhin2016}, thus illustrating the validity of Assertion~3~and~4.

As mentioned before, the evaluation of the inverse hysteresis operator requires the solution of a coupled minimization problem for all pinning forces, while the inner problems of the forward operator decouple. 
In Table~\ref{table:1}, we compare the computational efficiency of the two models 
in dependence of the number of pinning forces.
\begin{table}[ht!]
\caption{Computation times\label{table:1}}
\begin{center}
\begin{tabular}{c||c|c||c|c||c} 
 & \multicolumn{2}{c}{forward operator} & \multicolumn{2}{c}{inverse operator} & time quotient \\
$N_{\chi}$ & time $[ms]$ & it & time $[ms]$ & it & $\frac{\text{inverse}}{\text{forward}}$ \\[2pt]
\hline
\hline
$2$    & $0.253$ & $6.57$ & $0.724$ & $6.05$ & $2.99$    \\
$5$   & $0.386$ & $7.99$ & $1.74$  & $7.56$ & $4.66$ \\
$10$   & $0.420$ & $8.29$ & $4.10$ & $8.31$ & $9.84$ \\
$15$   & $0.483$ & $8.33$ & $7.94$ & $8.80$ & $16.3$ \\
$20$   & $0.523$ & $8.42$ & $12.4$ & $8.88$ & $25.1$ \\ 
\end{tabular}
\end{center}

Average computation times and iteration numbers (it) per load-step for hysteresis operators with varying number of pinning forces $N_{\chi}$
\end{table}
As expected, the computation times for the forward model increases linearly with the number $N_\chi$ of pinning forces. 
For the inverse hysteresis operator, a cubic dependence on $N_\chi$ can be expected. 
For both models, the number of Newton iterations for solving the minimization problems, on the other hand, seems not to be affected by the number of pinning forces.

\section{Magnetic field computation}

We now discuss in some detail the applicability of the  proposed forward and inverse hysteresis models for finite element simulations of scalar and vector potential formulations of a typical magnetic field problem. The description of the material models via energy resp. co-energy functionals allows to guarantee well-posedness of the resulting field equations and view them as convex minimization problems. As a consequence, globally convergent iterative solvers can be devised~\cite{Egger2024quasi,Egger2024vec}. 
These are also used in our computational tests below.

\subsection{General problem setup}
\begin{figure}[ht!]
    \centering
    \bigskip
    \includegraphics[trim={2cm 0.8cm 2cm 0.8cm},clip,width=0.4\textwidth]{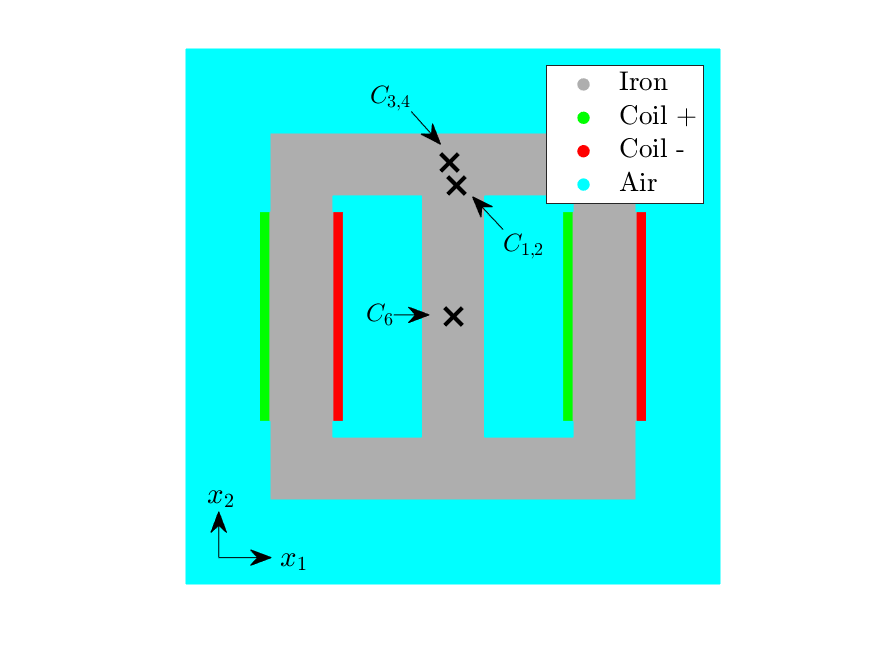}
    \caption{Sketch of the geometry of the TEAM~32 problem with iron (grey), coil (green and red), air (cyan). Evaluations of the magnetic field and flux at the points $C_6$, $C_{1,2}$ and $C_{3,4}$ are considered later on.}
    \label{fig:geometry}
\end{figure}

We consider the TEAM~32 benchmark problem~\cite{team32}. The relevant two-dimensional cross-section of the computational domain is depicted in Fig.~\ref{fig:geometry}. 
A typical load cycle is considered, in which the driving current in the coils is altered periodically.
Following~\cite{Meunier2008}, the quasi-static fields in the $i$-th time step $t^i$ are assumed to satisfy
\begin{alignat}{2}
\curl \H^i = j_s^i \quad \text{in } \Omega&, \qquad & \div \B^i &= 0 \quad \text{in } \Omega,\label{eq:magneto} \\
\B^i = \grad_\H w_*^i(\H^i)&, \quad & \B^i \cdot \n &= 0 \quad \text{on } \partial\Omega,\label{eq:magneto2}
\end{alignat}
with co-energy density $w_*^i(\H)$ described below. 
We will also make use of the alternative representation
\begin{alignat}{2}
\H^i = \grad_\B w^i(\B^i) \label{eq:magneto2a}
\end{alignat}
of the material law which is equivalent due to \eqref{eq:8}
In the two-dimensional setting under consideration, only the in-plane components of $\H$ and $\B$ are considered, and $j_s^i$ amounts to the out-of-plane component of the current.
The current density in the coils is given by $j_s^i=j_0 \cdot I_s(t^i)$ with $j_0 = \pm \frac{1}{A_w}$, $A_w = 5\times10^{-4}\,m^2$ having unit average but different signs on the two coil areas. Hence $I_s$ amounts to the total current flowing through the coils with  $N_w = 90$ number of turns.
A piecewise linear source field $\H_0$ with $\curl \H_0 = j_0$ and finite support in a vicinity of the coils is computed analytically~\cite{Engertsberger23} and used throughout our computations.  

\subsection{Material models}
For air and coil regions, we use a linear material model with (co)-energy densities $w^i(\B)=\frac{1}{2\mu_0} |\B|^2$ resp. $w^i_*(\H) = \frac{\mu_0}{2}|\H|^2$ independent of the time step.
The ferromagnetic behavior of iron, on the other hand, is described by the equivalent forward or inverse hysteresis models 
$w_*^i(\H)=w_*(\H,\{\J_p^i\})$ or 
$w^i(\B)=w(\B;\{\J_p^i\})$
as defined in \eqref{eq:13} and \eqref{eq:15}.  
In our computations, we use a model with $5$ pinning forces and parameters $\chi_k, w_k$ chosen as in~\cite{Lavet2013}.
The internal energy densities are defined by~\eqref{eq:Uk} with $A = 90.302$ and $J_s = 1.573$\,T. 
The magnetic polarizations are initialized by $\J_{k,p}^{0}=0$. 

\subsection{Scalar potential formulation}
The magnetic field is split by $\H^i = \H_s^i - \nabla \psi^i$ into a source field $\H_s^i = \H_0 \cdot I_s(t^i)$ and a reaction field which is represented by a reduced scalar potential $\psi^i$. 
Using this ansatz \eqref{eq:magneto}--\eqref{eq:magneto2} leads to the simplified system 
\begin{alignat}{2}
    \div(\grad_\H w^i_*(\H_s^i - \nabla \psi^i) &= 0 \quad &&\text{in } \Omega, \label{eq:scalar1}\\
    \n \cdot (\grad_\H w^i_*(\H_s^i - \nabla \psi^i) &= 0 \quad && \text{on } \partial\Omega, \label{eq:scalar2}
\end{alignat}
which naturally involves the \emph{forward} hysteresis operator~\eqref{eq:14} with co-energy density defined by~\eqref{eq:13}.
These equations are the first order optimality conditions for the convex minimization problem
\begin{align}
\min_{\psi} \int_\Omega w^i_*(\H_s^i - \nabla \psi) \, dx.
\end{align}
Standard $H^1$-conforming finite elements can be used for discretization of $\psi^i$, leading to a large scale optimization problem resp. nonlinear system to be solved in every time step. 
Let us note that the nonlinear term in \eqref{eq:scalar1}--\eqref{eq:scalar2} is not differentiable, so the classical Newton method cannot be applied here.
Globally convergent iterative solver based on local Quasi-Newton updates can, however, be devised for the efficient numerical solution; see \cite{Egger2024quasi} for a detailed analysis.

\subsection{Vector potential formulation}
The magnetic flux is represented by $\B^i = \Curl a^i$, where $\Curl$ is the scalar-to-vector operator and $a^i$ denotes the out-of-plane component of the magnetic vector potential. Inserting this ansatz into \eqref{eq:magneto}--\eqref{eq:magneto2a} leads to
\begin{alignat}{2}
    \curl(\grad_\B w^i(\Curl a^i) &= j_s^i \quad &&\text{in } \Omega, \\
    a^i &= 0 \quad &&\text{ on } \partial \Omega,
\end{alignat}
now involving the \emph{inverse} hysteresis operator~\eqref{eq:15}. 
These equations can again be interpreted as the first order optimality system of a convex minimization problem
\begin{align}
\min_a \int_\Omega w(\Curl a) - j_s^i \; a \, dx,
\end{align}
which has to be solved in every time step.
Standard $H^1$-conforming finite elements can be used for the discretization of $a^i$, and the resulting large scale optimization problems resp. nonlinear systems are solved iteratively with the same methodology as employed for the scalar potential formulation. Global convergence of the resulting schemes can be established~\cite{Egger2024quasi,Egger2024vec}. 

\subsection{Numerical results}
\begin{figure}[ht!]
    \centering
    \includegraphics[trim={0cm 0.2cm 0cm 0.6cm},clip,width=0.48\textwidth]{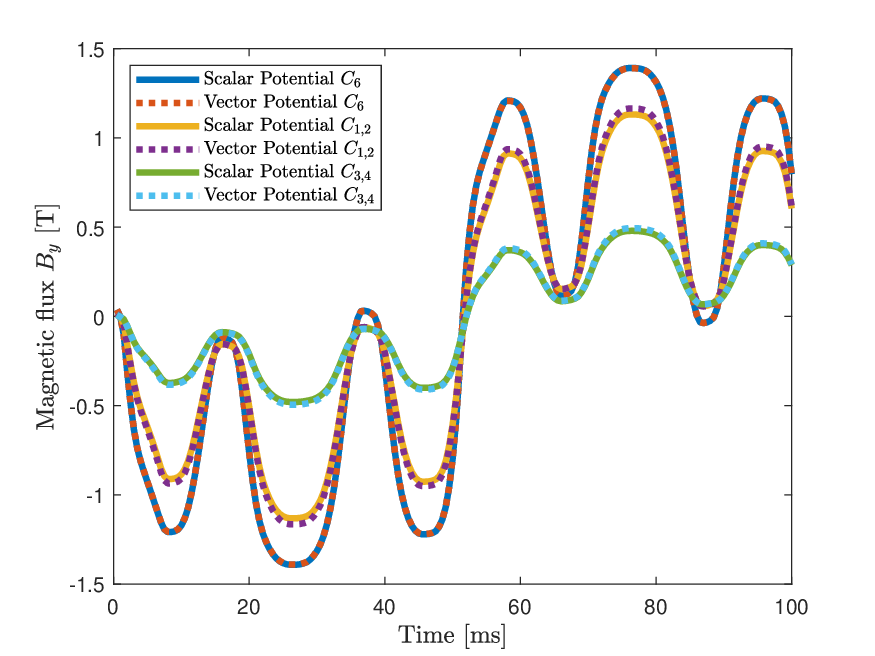}
    \caption{Vertical component $B_y(t)$ of magnetic flux at the points $C_6$, $C_{1,2}$ and $C_{3,4}$ for the scalar potential and vector potential formulation.
}
    \label{fig:magnetic:fluxes}
\end{figure}

For our numerical tests, we consider the test case~2 of TEAM~problem~32, see \cite[Fig.~7]{team32}, which prescribes a sinusoidal current excitation $I_s(t)$ with some additional higher harmonics.
In Fig.~\ref{fig:magnetic:fluxes}, we display the vertical component $B_y^i\simeq B_y(t^i)$ of the magnetic flux at selected points $C_6$, $C_{1,2}$, and $C_{3,4}$ of the geometry, obtained by finite element simulations for the scalar and vector potential formulations, respectively.
Despite the fact that different formulations and finite element approximations are used, the simulations based on the scalar and vector potential yield practically identical results.
The small difference can in fact be attributed to discretization errors which vanish when further refining the mesh. On the continuous level, the two formulations are equivalent.

\begin{figure}[ht!]
    \centering
    \includegraphics[trim={0cm 0cm 0cm 0.2cm},clip,width=0.48\textwidth]{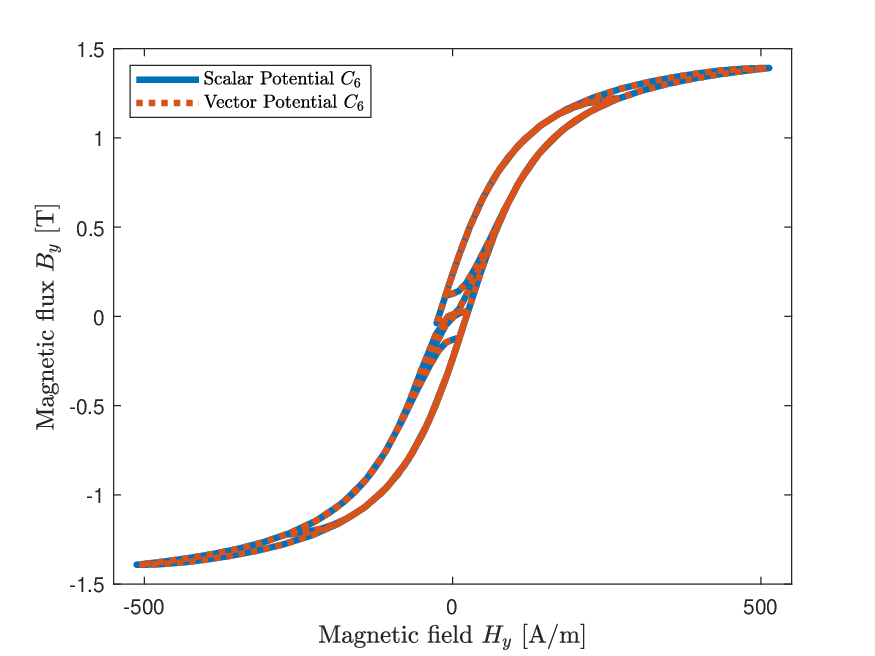}
    \caption{Hysteresis loop at the geometry point $C_6$ computed by the scalar potential and vector potential formulation, respectively.}
    \label{fig:6}
\end{figure}

For completeness of the presentation, we further depict in Fig.~\ref{fig:6} the hysteresis loops observed at the point $C_6$ of the geometry, which again illustrates the equivalence of the underlying forward and inverse hysteresis operators. 

\section{Discussion}

The magnetic vector hysteresis model of Bergvist~\cite{Bergqvist1997,Bergqvist1997a}, further developed by Henrotte et al.~\cite{Henrotte2006,Lavet2013}, implicitly describes magnetic flux $\B=\B(\H,\J_p)$ as a function of magnetic field $\H$ and previous magnetic polarization $\J_p$. By defining an 
appopriate magnetic co-energy density $w_*(\H;\J_p)$, a reformulation of this model in the form $\B=\grad_\H w_*(\H;\J_p)$ was derived. 
Using arguments of convex duality, a compact form of the corresponding energy density $w(\B;\J_p)$ could be given which allowed us to explicitly describe the inverse hysteresis operator in the form $\H=\partial_\B w(\B;\J_p)$.
The well-posedness of the two material models, their equivalence, and efficient implementation has been discussed theoretically and demonstrated numerically.
In addition, the application of the forward and inverse hysteresis operators in finite element approximations of scalar or vector potential formulations has been demonstrated by simulations for a typical benchmark problem.
The representation of the material models via appropriate energy respectively co-energy densities was the key ingredient for the analysis and efficient numerical realization.

%
The evaluation of the forward and inverse hysteresis operators requires the multiple solution of non-smooth minimization problems at every quadrature resp. material point. In the case of single pinning forces, this can be done efficiently using the approach discussed in~\cite{Prigozhin2016}. If the minimization problems for different pinning strengths separate, this can also be extended to the case of multiple pinning forces.
The efficient solution of the coupled non-smooth minimization problems resulting otherwise, however, deserves further investigation. 
A convergence analysis of iterative solvers for the nonlinear systems arising after finite element approximations of scalar or vector potential formulations using the proposed hysteresis operators has been presented in~\cite{Egger2024quasi}. 
The mathematical justification of some of the technical assumptions used in this analysis and the further acceleration of the nonlinear solvers by use of semi-smooth Newton methods are left as topics for future research.

\begingroup
\section*{Acknowledgements}
This work was supported by the joint DFG/FWF Collaborative Research Centre CREATOR (DFG: Project-ID 492661287/TRR 361; FWF: 10.55776/F90).
\endgroup


\end{document}